\documentclass[12pt]{article}
\usepackage{epsfig}
\usepackage{amsmath}
\usepackage{amsthm, amssymb, amscd}

\newcommand{\vn}{$\{v_1, \dots v_n \}$}
\newcommand{\starv}{$\{\nu_1, \dots \nu_n, \nu_{n+1} \}$}
\newcommand{\stare}{$\{\epsilon_1, \dots \epsilon_n \}$}
\newcommand{\ej}{$\{e_{\{i,j\}}, 1 \leq i < j \leq n \}$}

\newcommand{\Z}{{\hspace{-0.3pt}\mathbb Z}\hspace{0.3pt}}
\newcommand{\R}{{\hspace{-0.3pt}\mathbb R}\hspace{0.3pt}}

\author{Hugh Nelson Howards \\ Wake Forest University; and \\
Visiting Researcher \\ University of Texas, Austin\\
e-mail: howards@mthcsc.wfu.edu}
\date{}

\title{Knotted Spheres and Graphs in Balls}

\begin{document}
\maketitle

\newtheorem{definition}{Definition}[section]
\newtheorem{corollary}[definition]{Corollary}
\newtheorem{corrolary}[definition]{Corollary}
\newtheorem{theorem}[definition]{Theorem}
\newtheorem{proposition}[definition]{Proposition}
\newtheorem{lemma}[definition]{Lemma}
\newtheorem{claim}[definition]{Claim}
\newtheorem{ex}[definition]{Example}
\newtheorem{q}[definition]{Question}
\newtheorem{exer}[definition]{Exercise}

\begin{abstract}
Given a properly embedded graph $\Gamma$ in a ball $B$ and a 
punctured sphere $\Sigma$ properly embedded in $B- \Gamma$
we examine the conditions on $\Gamma$ that are necessary to 
assure that $\Sigma$ is boundary parallel.

\medskip

\medskip
 

\end{abstract}


\section{Knotted Spheres}
\label{section:knottedspheres}
This paper explores embedding punctured spheres in balls.
I would like to thank Mike Freedman for suggesting the question
and for insightful comments.
I would also like to thank Cameron Gordon,
John Luecke, Martin
Scharlemann, and Ying-Qing Wu for helpful comments along the way.

We first must
 set up some definitions that greatly simplify the
statements of the theorems. Throughout the paper
let $B$ be the unit ball in $\R^3$.
Let $S$ be the boundary of $B$.

\begin{definition} A complete graph on $n$ vertices
whose vertices \vn
are disjoint disks on $S$, and whose edges \ej are 
properly embedded geodesics
in $B$ with
the property that if $s \neq t$, $v_s \cap v_t = \emptyset$, 
is a \underline{standard unlinked n-graph}.
\end{definition}

\begin{definition} A complete graph on $n$ vertices
whose vertices \vn are disjoint disks on $S$ and whose edges 
\ej are 
properly-embedded, disjoint unknotted arcs
in $B$ with
the property that no two edges are linked in the ball, 
is called a \underline{pairwise unlinked n-graph}.
\end{definition}

\begin{definition}  Given a graph $\Gamma$ in $B$,
and $\Sigma$ a properly embedded  n-holed sphere in $B$,
such that $\partial \Sigma = \cup_{i=1}^n \partial v_i$ and 
$\Sigma \cap e_{\{i,j\}} = \emptyset$ for all $i$, $j$.
$\Sigma$ is called an \underline{enveloping n-holed sphere}
for $\Gamma$.  $B-\Sigma$ consists of two components. The part containing
\ej is called the inside of $\Sigma$, the other component is 
called the outside.

\end{definition}

\begin{definition}
An enveloping n-holed sphere $\Sigma$ is \underline{standard}, if 
it is boundary parallel, meaning there is 
a product structure on the outside of $\Sigma$ taking the interior of
$\Sigma$ to $S-\cup_{i=1}^n v_i$, but leaving $\partial \Sigma$ fixed.
\end{definition}

\begin{definition}
Let a \underline{core} of an enveloping n-holed sphere 
$\Sigma$ be a 1-complex,
such that the boundary of a regular neighborhood of the
complex co-bounds a product region
with $\Sigma$.
Let the \underline{star core} be the
core with exactly $n$ edges and one vertex of valence 
of $n$. We shall call the 
vertices of the star core \starv where $\nu_i \subset v_i$ and $\nu_{n+1}$
is the vertex of valence $n$.  We call the
edges of the star core \stare where  $\epsilon_i$ is the edge 
containing  $\nu_i$.

\end{definition}

\section{Pairwise unlinked n-graphs}
\subsection{The Core Lemmas}

In this section we introduce a couple of basic lemmas that will 
increase our insight to the theorem and will be useful
in some of the cases.

\begin{lemma}
\label{lemma:knottedcore1}
Given $\epsilon_i$ and $\epsilon_j, i \neq j$,
two edges of the star core of an enveloping
n-holed sphere for a pairwise unlinked n-graph,
$\epsilon_i \cup \epsilon_j$ may never be a knotted arc in $B$. 
\end{lemma}

\begin{figure}[htbp] 
\centerline{\hbox{\psfig{figure=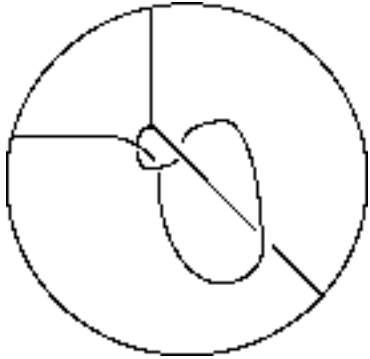}}}
\caption{A knotted core}\label{fig:knottedcore}
\end{figure}

\begin{proof}  $\epsilon_i \cup \epsilon_j$
makes up the core of the cylinder
that is left over if the n-holed sphere is compressed
along disks parallel to all of the uninvolved vertices.
If the core of the cylinder were knotted (see 
Figure~\ref{fig:knottedcore}),
then by standard
satellite knot theory so is any arc running through the
cylinder, but the edge between these two vertices in our
graph is unknotted and can be assumed to be
inside the cylinder, so this cannot be the case.
\end{proof}

\begin{lemma}
\label{lemma:linkedcore1}
Let $\epsilon_i$, $\epsilon_j$, and $\epsilon_k$ be three
edges of a star core of an enveloping n-holed sphere
for a pairwise unlinked n-graph.  Let $\epsilon_i^k$ and $\epsilon_j^k$
be edges obtained from $\epsilon_i$ and $\epsilon_j$
by contracting $\epsilon_k$ and perturbing the arcs slightly so their end
points are disjoint.  Then $\epsilon_i^k$ and $\epsilon_j^k$
are not linked.

\end{lemma}
\begin{proof} If $\epsilon_i^k$ and $\epsilon_j^k$
were linked, then $e_{\{i,k\}}$ and $e_{\{j,k\}}$
(edges of the 
the original unlinked n-graph running from vertex $i$ to vertex $k$
and $j$ to $k$ respectively)
would have to be linked.  

\end{proof}

[C] and [GF] are good places to look for an introduction to 
rational tangles. 
As explained in [GF] a rational tangle $T$ is assigned a 
rational number $F(T)$
corresponding to a simple continued fraction.  
Two tangles $T_1$, $T_2$ are isotopic if and only if 
$F(T_1) = F(T_2)$.
Let $N(T)$ be the knot obtained by from $T$ by connecting the 
ends of the tangle 
in the manner called the numerator of $T$ and $D(T)$ the denominator 
(See Figure~\ref{fig:tangles}).

\begin{lemma}
If $T$ is a rational tangle with $D(T)$ an unknot, then
$F(T) = p, p \in \Z$.  Therefore $T$
may be undone by merely twisting two of the vertices around
each other $p$ times leaving the other two vertices fixed.
\label{lemma:rationaltangles1}
\end{lemma}
\begin{proof}  

It is well known that if $F(T)=p/q$, then $D(T)$ 
produces a $p/q$ 2-bridge link, which is trivial if and only if q=1.
Similarly, of course, $N(T)$ is
an unknot if and only if $p=1$ so in that case $F(T)= 1/q, q \in \Z$ (the
picture is just rotated ninety degrees).
[M] section 9.3 and
[BZ] sections 12.A and 12.B
have expositions on the these facts.
The classification was first done in [S].

\begin{figure}[htbp] 
\centerline{\hbox{\psfig{figure=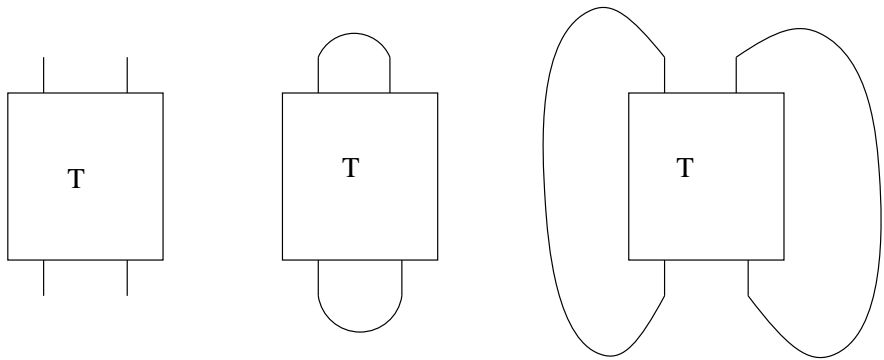}}}
\caption{A representation of tangle $T$ is given on the left.  
The numerator of $T$ is demonstrated in the center,
the denominator on the right.}\label{fig:tangles}
\end{figure}
\end{proof}

\subsection{The case, $n=3$.}

\begin{figure}[htbp] 
\centerline{\hbox{\psfig{figure=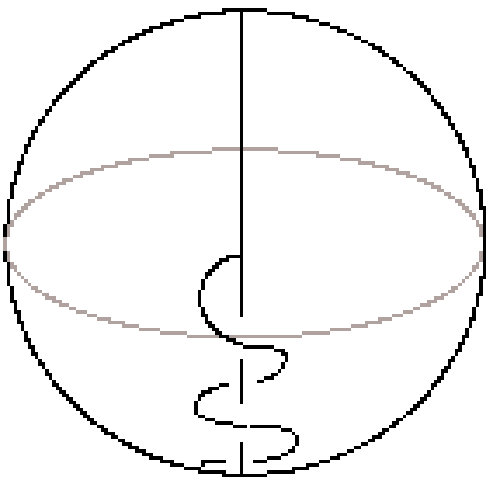}}}
\caption{The only possible complication for $n=3$.}\label{fig:gen3star}
\end{figure}

\begin{theorem}
Every enveloping 3-holed sphere for a pairwise unlinked 3-graph
is standard. \label{thm:threearcs}
\end{theorem}
\begin{proof} By Lemma~\ref{lemma:knottedcore1} the star core is not knotted,
so we may assume it consists of one straight edge 
$\epsilon_1 \cup \epsilon_2$ running from
$\nu_1$ the north pole to $\nu_2$ at
the south pole and another edge $\epsilon_3$
meeting this
edge at $\nu_4$ at
the origin and winding around in some manner before
ending up at a $\nu_3$ somewhere on the 
Southern hemi-sphere (as in Figure~\ref{fig:pole}).

\begin{figure}[htbp] 
\centerline{\hbox{\psfig{figure=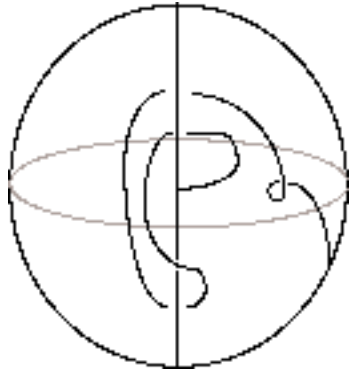}}}
\caption{Straightening the core}\label{fig:pole}
\end{figure}

Act upon $\epsilon_3$ by an ambient isotopy 
leaving $\epsilon_1 \cup \epsilon_2$ fixed
until $\epsilon_3$ lies 
entirely below the equatorial disk except for
at $\nu_4$
(Figure~\ref{fig:gravity}).
Now if we examine the two edges of the core
in $B'$ the Southern hemi-ball (the ball we get by
cutting the original ball in half along the 
equatorial disk) and pull  $\epsilon_2$ and  $\epsilon_3$
apart slightly so that they no longer intersect
to get  $\epsilon_2'$ and  $\epsilon_3'$, 
the core's edges cannot be linked by Lemma~\ref{lemma:linkedcore1}
(this is the same as contracting $\epsilon_1$ and looking at
$\epsilon_2^1$ and $\epsilon_3^1$.
This together with the fact that neither  $\epsilon_2$ nor
$\epsilon_3$ can be knotted in $B'$ by Lemma~\ref{lemma:knottedcore1},
means that  $\epsilon_2' \cup \epsilon_3'$
is just a rational tangle in $B'$.

\begin{figure}[htbp] 
\centerline{\hbox{\psfig{figure=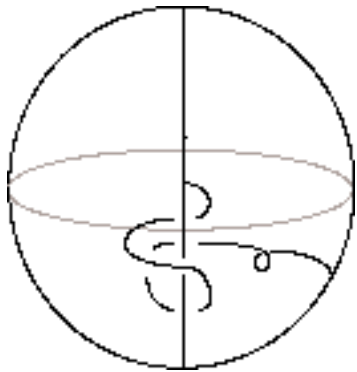}}}
\caption{An isotopy acting upon the straightened core}\label{fig:gravity}
\end{figure}

Now by Lemma~\ref{lemma:knottedcore1} we know 
that the rational tangle to which they
correspond,
must be an unknot if vertices on the equatorial disk are connected
as are the vertices  on the sphere.
By Lemma~\ref{lemma:rationaltangles1} we may assume 
that the rational tangle is obtained
merely by twisting the two vertices on the southern hemisphere around
each other, and therefore that it may be untangled without 
affecting the rest of the core (Figure~\ref{fig:gen3star}).
Therefore the core is standard, as must be the enveloping
n-holed sphere.
\end{proof}

\subsection{The Case $n = 4$}

The argument for the case $n = 4$ is particularly interesting, 
because the methods for the previous case fall slightly
shy of working, but the counterexample from the case $n \geq 5$
also just fails to disprove it.  I would
like to thank Ying-Qing Wu for ideas that were particularly
helpful in this section.

\begin{theorem}
Given $n \leq 4$ 
every enveloping n-holed sphere for a pairwise unlinked n-graph $\Gamma$
is standard. 
\label{thm:fourarcs}
\end{theorem}
\begin{proof} Assume there is a counterexample
for $n=4$ and examine what
it must look like.  $B$ minus the interior of the non-standard ball,
bounded by the non-standard enveloping 4-holed sphere
is homeomorphic to  $S^2 \times I$ minus 
an open neighborhood of four arcs 
$\{\epsilon_1, \epsilon_2, \epsilon_3, \epsilon_4 \}$ which
run from the inner sphere to the outer sphere.  Since 
the enveloping 4-holed sphere is non-standard, the arcs cannot
each simultaneously be isotoped in $S^2 \times I$
to be $pt \times I$.

Since the original edges of $\Gamma$ were 
unknotted and pairwise unlinked, 
they could be thought of as rational tangles.
Examine the three pairs of edges of $\Gamma$ corresponding 
to the three possible pairings of the vertices, 
$(e_{\{1,2\}}, e_{\{3,4\}})$,
$(e_{\{1,3\}}, e_{\{2,4\}})$, $(e_{\{1,4\}}, e_{\{2,3\}})$.  
Pick one pair, say $(e_{\{1,2\}}, e_{\{3,4\}})$.
Because the pair is a rational tangle, there is a 
disk $D_1$ in $B$ that separates
the edges. Assume $D_1$ intersects $\Sigma$ minimally.  
An innermost curve of intersection on $D_1$
yields a subdisk $D_1'$ which
compresses $\Sigma$ into two annuli $A_{\{1,2\}}$ and 
$A_{\{3,4\}}$ 
which are just the boundary of
a regular neighborhood of $e_{\{1,2\}}$ and $e_{\{3,4\}}$ respectively.  
The same arguments can be made for the other two pairs.
Thus, there are three disks that can be added to 
$S^2 \times 0$ that separate the end points of
the $\epsilon_i$ into pairs, each yielding a rational tangle.

If we take the branched double cover of $S^2 \times I$ over the 
four arcs $\epsilon_1, \epsilon_2, \epsilon_3, \epsilon_4$, 
we get a manifold $M$ with boundary two tori $T_0$, $T_1$
because the branched double cover of a sphere over four points 
is a torus.  
Since adding $D_1'$ to $B$ minus the
 ``inside'' of $\Sigma$, thought of as 
$S^2 \times I$ minus a neighborhood of the four branching arcs
yields the exterior of 
a rational tangle, and the branched double cover of a ball over a rational
tangle is, of course, just a solid torus. 
Thus, $D_1'$ lifts to
a disk that gives a filling of $T_0$ that turns $M$ into a solid torus,
so $M$ is just a solid torus minus a knot.
Note that the three pairings of
the vertices gives three different fillings of $T_0$ 
each of which yields a solid torus.

\begin{theorem} {\bf [B Corollary 2.9]}
If $k$ is a nontrivial knot in $D^2 \times S^1$ such that $k$
is not parallel into $\partial D^2 \times S^1$ and 
there exists more than one nontrivial surgery on 
$k$ yielding $D^2 \times S^1$, then $k$ is equivalent either to
$W^1_3W^{-3}_7$ or its mirror image $W^{-1}_3W^{3}_7$.
\end{theorem}

$W^{-1}_3W^{3}_7$ is the $(-2,3,7)$ pretzel knot embedded as
shown in Figure~\ref{fig:berge}. We refer to the this as the 
Berge knot. Since the arcs were not standard, we must lift to the
Berge knot or a knot $k$ parallel into $T_1$, the boundary of 
$D^2 \times S^1$.

\begin{proposition}
If $k$ is parallel into the boundary of $D^2 \times S^1$ then $k$ fails
to produce a counter example to Theorem~\ref{thm:fourarcs}.
\end{proposition}

The proof of the proposition requires two steps.  First we prove that 
$k$ can be assumed to have the standard embedding in $D^2 \times S^1$
by showing there is a unique strong inversion on $k$.  Second we prove
that the punctured sphere produced by quotienting out by the $\Z_2$ symmetry
is either boundary parallel or else
violates Lemma~\ref{lemma:knottedcore1} and therefore is not a counterexample.

\begin{proof}

To prove that $k$ must be embedded
in the standard manner we examine an annulus $A$ that runs from
$T_0$, the torus boundary component corresponding to $k$ 
in the exterior of $k$, to $T_1$, the boundary of the solid torus.
$A$ can be assumed to be embedded (see [CF]). Let $F$ be the 
involution of $D^2 \times S^1$. Let $F(A)=A'$. 

\begin{lemma}
There is a unique strong inversion of $k$, a torus knot parallel
to the boundary of a solid torus.
\end{lemma}
\begin{proof}
Our first goal is to show that 
$A$ can be chosen such that $A=A'$.  We choose $A$ to have a minimal number
of intersections with $A'$. It is clear that $\partial A$
can be assumed to be fixed by $F$, so perturbing $A$ slightly we can assume
that $\partial A \cap \partial A' = \emptyset$.  Now $A \cap A'$ must consist
of simple closed curves.  An innermost disk argument suffices to show that 
each of the simple closed curves is essential on $A$ and $A'$,
so the intersection consists of 
disjoint simple closed curves $c_1, c_2, \dots c_n$ parallel to $\partial A$
on $A$.  Let $c_1$ be the curve of intersection closest to 
the boundary component of $A$ on $T_1$, $c_2$ be the second, and so on 
increasing the index as the curves move towards the boundary component
on $T_0$. Likewise on $A'$ the intersection consists of 
parallel essential 
circles $c_1', c_2', \dots c_n'$ labeled in the same manner. 
Let $A_i$ be the sub-annulus of $A$ running from the boundary 
component on $T_1$
to $c_i$, and $A_i'$ be the sub-annulus of $A'$ running 
from the boundary component on $T_1$ to $c_i'$.  Let $c_j$ be 
the curve of intersection 
on $A$ that corresponds to the intersection 
with $A'$ at $c_1'$. If we cut and paste
$A$ replacing $A_j$ by a push off of $A_1'$ 
we reduce the number of intersections
of $A$ and $A'$ by at least one.  In order to 
preserve the property that $F(A)=A'$
we must also replace $A_j'$ by a push off of 
$A_1$.  This, however, cannot increase the 
number of intersections, so we have a new 
annulus running from $T_1$ to $T_0$ that has
fewer intersections with its image under 
$F$, contradicting minimality.  Thus,
we can assume that $A \cap A' = \emptyset$.

This, however, implies
that restricting to the solid torus between $A$ and $A'$, 
$F$ takes $A \times I$ to 
itself, exchanging $A \times 0$ with $A \times 1$.  This in 
turn implies that there 
must be an annulus in $A \times I$ that is fixed by $F$. 

Now that we know that $A$ is fixed, we use it to show that we 
have a standard inversion
of $k$.  Let $D$ be a meridional disk for $T_1$ that is fixed 
by $F$.  Examine $D \cap A$.
In general the intersection pattern on $D$ will look something like 
Figure~\ref{fig:diskannulus} with a collection of arcs running 
from $T_0$ (which punctures $D$ several times) to $T_1$
and a collection that run from 
one of the punctures from $T_0$ to another. We can remove the arcs 
running from 
$T_0$ to $T_0$ by picking 
an outermost arc on $A$ that runs from one component of 
$\partial A$ to itself. This small disk gives an isotopy of $A$ together with
$k$ that reduces the number of intersections of $A$ with $D$.  
because $F(A)=A$
we can simultaneously do a second isotopy of $A$ that also reduces the number 
of intersections of $A$ with $D$ and preserves the symmetry of $k$.  Thus,
we can assume that $D \cap A$ consists solely of arcs running 
from $k$ to $T_1$.
This, however, shows that we have the standard symmetry for $k$ 
because cutting
the solid torus along $D$  turns it into a cylinder and 
$A$ becomes
bands running from the top of the cylinder to the bottom in the 
unique way possible.
\end{proof}

\begin{figure}[htbp] 
\centerline{\hbox{\psfig{figure=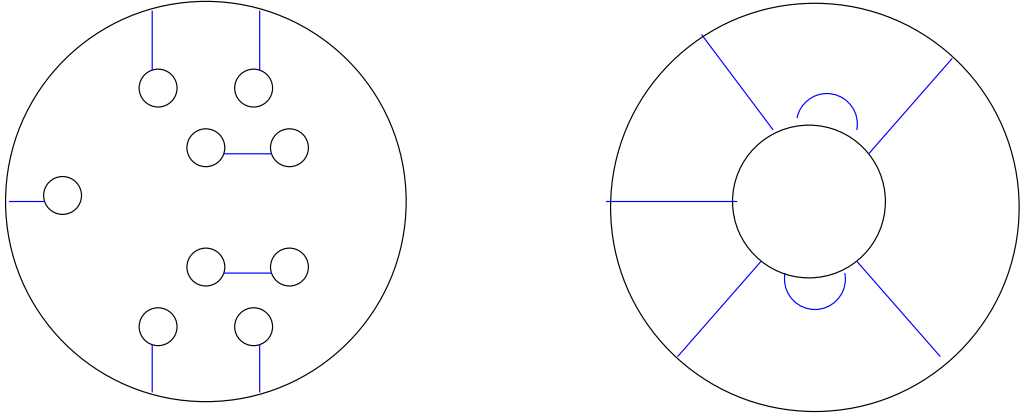}}}
\caption{$D \cap A$. $D$ is pictured punctured 
nine times by $T_0$ 
on the left and $A$ is pictured on the right.}
\label{fig:diskannulus}
\end{figure}

Now we need only argue that torus knots with standard embeddings 
fail to give a counterexample.  
Let $k$ be a  torus knot embedded in $D^2 \times S^1$,
fixed by an involution $F$ of $D^2 \times S^1$.  Let $A$ be 
the fixed annulus above.
Let $M$ be $S^2 \times I$ with four branching arcs
$\epsilon_1, \dots \epsilon_4$, the quotient of 
the exterior of $k$ in $D^2 \times S^1$
by $F$.  In the quotient, $T_0$ maps down to 
$S^2 \times 0$ which we will designate $S_0$
and $T_1$ maps to $S^2 \times 1$ designated $S_1$.

\begin{lemma}
Either $S_0 \cup \epsilon_1 \cup \epsilon_2 \cup \epsilon_3
\cup \epsilon_4$
violates Lemma~\ref{lemma:knottedcore1} 
or it is standard
and therefore in either case is 
not a counterexample to Theorem~\ref{thm:fourarcs}.
\label{lemma:nocounter}
\end{lemma}

\begin{proof}

Because 
$\epsilon_1, \dots \epsilon_4$, each run from $S_0$ to
$S_1$ 
the two arcs $a_1$, $a_2$ that are fixed in $D^2 \times S^1$ by
$F$ must each intersect $k$ in exactly one point.  In the knot
exterior $a_1$ and $a_2$ are then broken each into two arcs with one
end point of each arc on $T_1$ and one end point on $T_0$.
In $M$ the $a_i$ then become the four
branching arcs  $\epsilon_1 \dots \epsilon_4$.

Under the quotient, the annulus $A$ becomes a rectangular
disk $D$ which without loss of generality
runs down $\epsilon_3$
along $S_0$ up $\epsilon_4$ and back along $S_1$ and is disjoint
from $\epsilon_1$ and $\epsilon_2$.  This is clear because
$A \cap (a_1 \cup a_2)$ consisted of the
preimage of $\epsilon_3$ 
and $\epsilon_4$ but was disjoint
from the preimage of $\epsilon_1$ and $\epsilon_2$
(recall that $A$, like the $a_i$ is fixed by $F$ and runs
from $T_0$ to $T_1$).  
This means that $M$ looks exactly like Figure~\ref{fig:torustangle},
where $\epsilon_1$ and $\epsilon_2$ form a rational tangle
inside the ball designated $T$.

\begin{figure}[htbp] 
\centerline{\hbox{\psfig{figure=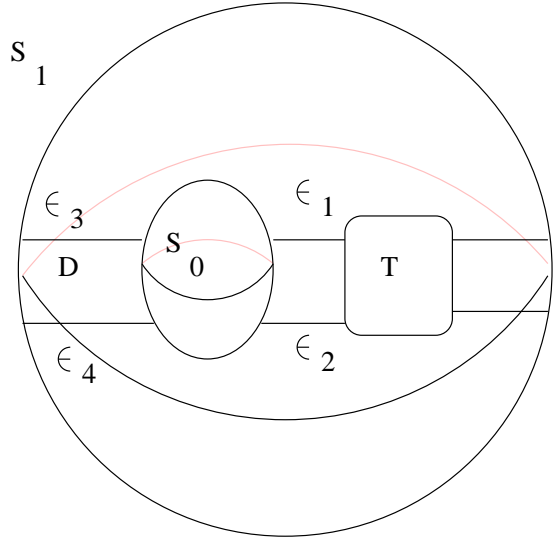}}}
\caption{The branching arcs of $M$ consist of two standard arcs 
$\epsilon_3$ and $\epsilon_4$ and a rational
tangle running from $S_0$ to $S_1$ designated $T$ consisting of
arcs $\epsilon_1$ and $\epsilon_2$.}
\label{fig:torustangle}
\end{figure}
 
By Lemma~\ref{lemma:rationaltangles1}
we see that $\epsilon_1 \cup \epsilon_2$ must be the tangle $T$ that 
results from two 
horizontal arcs whose eastern vertices are twisted
$n$ times around each other or $\epsilon_1 \cup S_0 \cup \epsilon_2$ will
be knotted violating Lemma~\ref{lemma:knottedcore1}.  
On the other hand, if $\epsilon_1 \cup \epsilon_2$ is the tangle
$T$ above, then twisting the eastern portion of $M$, $n$ times shows 
$M$ is homeomorphic to the standard picture and 
therefore again fails to be a counterexample.  

\end{proof}

{\em Note:} \hspace{.5mm} One can in fact prove that 
(if $k$ is a torus knot other than the unknot) $k$
never produces a standard enveloping sphere, but instead
always creates one that violates Lemma~\ref{lemma:knottedcore1}. 

\medskip

Lemma~\ref{lemma:nocounter} 
completes the proof that $k$ cannot be parallel to $T_1$ and 
therefore must be a Berge
knot if it is to produce a counterexample.
\end{proof}

\medskip

We may transform the traditional picture of the
Berge Knot to a symmetric one as in Figure~\ref{fig:berge}.
Snappea [W] tells us that this knot (entered
as a link) has exactly one
 $\Z_2$ symmetry, which we can now see.

\begin{figure}[htbp] 
\centerline{$\begin{array}{cc}
\psfig{figure=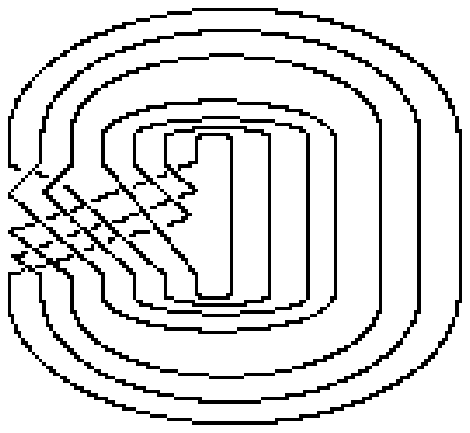}&
\psfig{figure=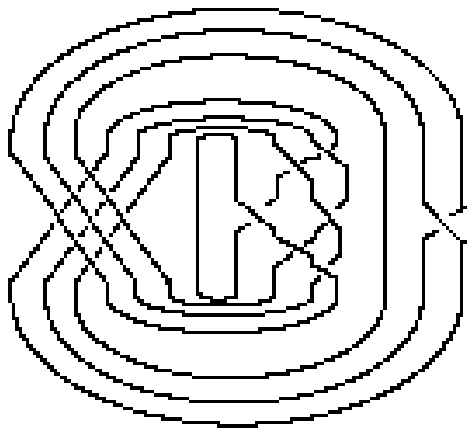}
\end{array}$}
\caption{The Berge knot
in $D^2 \times I$ and an embedding after 
an isotopy in $D^2 \times I$ that
reflects the knots $\Z_2$
symmetry ($D^2 \times I$ is not drawn, but is 
the obvious choice for the initial braid).}
\label{fig:berge}
\end{figure}

We quotient out by the symmetry to either get a counterexample, 
or proof that there is none.
We get $S^2 \times I$ minus four arcs
$\epsilon_1, \epsilon_2, \epsilon_3, \epsilon_4$. 
We will show that three
of the arcs violate Lemma~\ref{lemma:linkedcore1}
(See Figure~\ref{fig:linked1}).

\begin{figure}[htbp] 
\centerline{\hbox{\psfig{figure=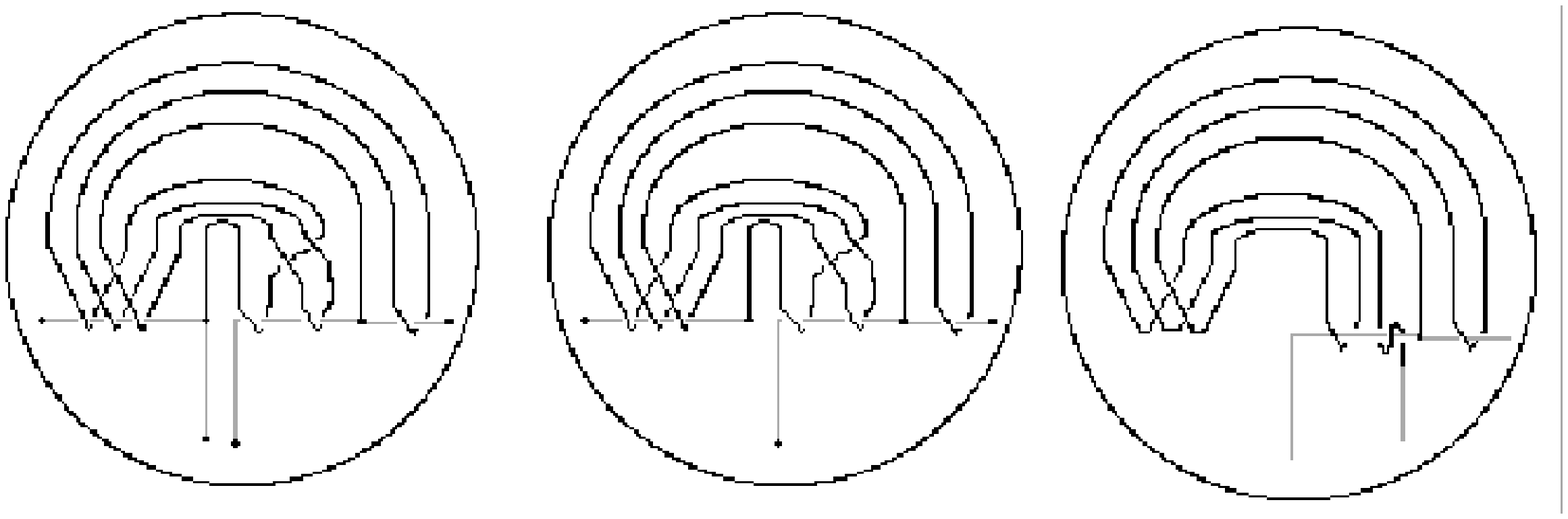,width=6in}}}
\centerline{\hbox{\psfig{figure=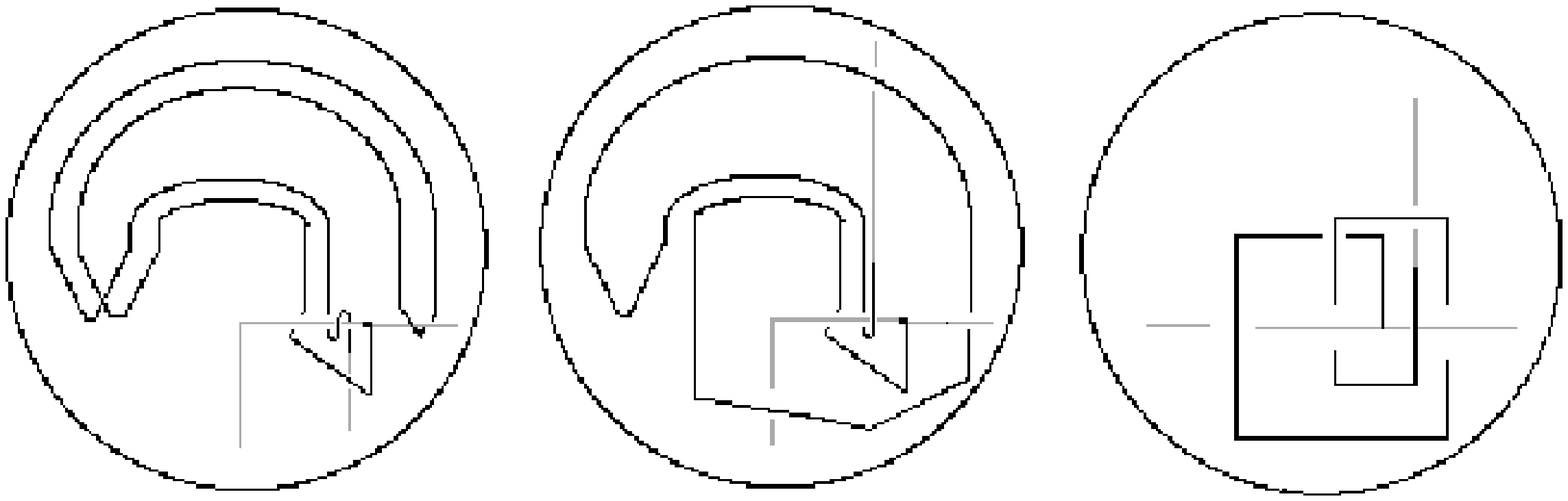,width=6in}}}
\caption{The dark edge is the interior sphere in $S^2 \times I$
obtained from quotienting out by the $\Z_2$ symmetry of the 
Berge knot in $D^2 \times I$.  The light edges are the four singular
arcs $\{a_1, a_2, a_3, a_4 \}$.  After dropping one of the 
$a_i$ and an isotopy, 
the three other edges appear to be linked.}
\label{fig:linked1}
\end{figure}

Label the edge omitted from the picture  $\epsilon_4$.
Label the horizontal edge that has one
of its vertices on the eastern side of the sphere
$\epsilon_3$.
Contraction of $\epsilon_3$ takes edges  $\epsilon_1$
and $\epsilon_2$ to $\epsilon_1^3$ and $\epsilon_2^3$
as pictured in Figure~\ref{fig:linked2}.

To see that $\epsilon_1^3$ and $\epsilon_2^3$
are in fact linked we call on 
the following well known facts about rational tangles. 
(See, for example,  [M]  {\bf Theorem 9.3.1}.)
\begin{theorem}
\begin{enumerate}
\item A 2-bridge knot (or link) is the denominator of some rational
tangle
\item Conversely, the denominator of a rational tangle is a 2-bridge knot
(or link).
\end{enumerate}
\label{thm:2bridge}
\end{theorem}

\begin{corollary}
If $k$ is the denominator
of a rational tangle, then $k$ is prime. 
\label{cor:primetangles}
\end{corollary}

The denominator of the tangle in Figure~\ref{fig:linked2} is
the connect sum of a trefoil and a figure eight knot and therefore
it is not a rational tangle by Corollary~\ref{cor:primetangles}.  
Since neither of the arcs is knotted, they must be linked
violating Lemma~\ref{lemma:linkedcore1}.
Since this was the only possible
counterexample, Theorem~\ref{thm:fourarcs} must be true. 

\begin{figure}[htbp] 
\centerline{\hbox{\psfig{figure=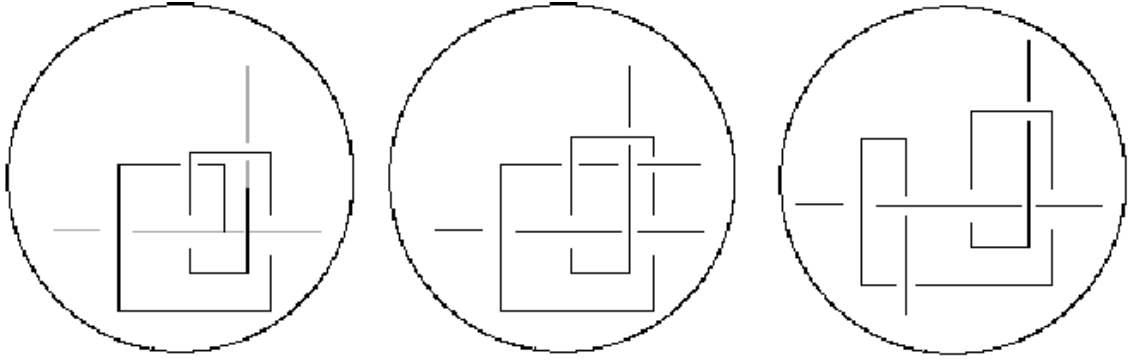,width=6in}}}
\caption{The edges violate Lemma~\ref{lemma:linkedcore1}.}
\label{fig:linked2}
\end{figure}
  
\end{proof}

\subsection{The Case $n \geq 5$}

\begin{theorem}
There exist pairwise unlinked n-graphs with non-standard
enveloping n-holed spheres for all $n$, $n \geq 5$. 
\label{fivearcs}
\end{theorem}
\begin{proof} 
Figure~\ref{fig:thetan} shows a $\theta_n$ curve
in a ball.  Such a graph is well known  not to be standard,
but every subgraph is standard.
Let a $\theta_n$ curve be the star core of 
an enveloping n-holed sphere. Although it is not standard,
it supports a pairwise unlinked n-graph $\Gamma$.
If we think of the enveloping sphere as bounding a 
central (round) ball with n tentacles running to
the boundary, we can picture the pairwise unlinked graph
as being a standard unlinked graph in the central ball
which is extended by a product down each of the tentacles.
Now since $n \geq 5$ any two edges of $\Gamma$ miss at least one
vertex, but there is an isotopy of any $n-1$ arcs
of a $\theta_n$ curve that makes those arcs appear
standard.  Likewise since the edges of $\Gamma$   completely
miss one of the tentacles, they may be pictured as being embedded
in a ball bounded by a standard enveloping $(n-1)$-holed sphere.
The arcs remain standard within the central ball and are
extended by a product down the tentacles throughout the entire
process, so clearly the edges are not linked pairwise.
(See Figure~\ref{fig:thetan}).

\begin{figure}[htbp] 
\centerline{\hbox{\psfig{figure=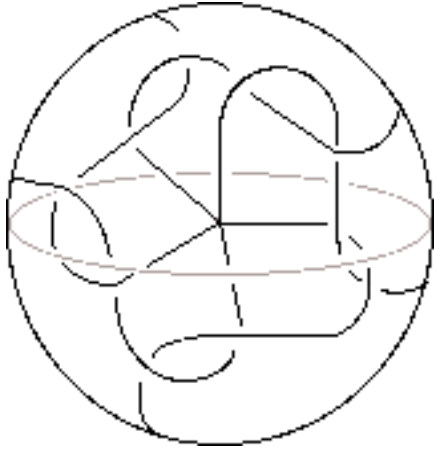}}}
\caption{The $\theta_n$ ($n= 5$ here)
curve is a non-standard core
but supports a pairwise unlinked graph for $n \geq 5$}\label{fig:thetan}
\end{figure}
\end{proof}

\section{Standard unlinked n-graphs}

\begin{theorem}
Given a standard unlinked n-graph $\Gamma$, 
every enveloping n-punctured sphere $\Sigma$  is standard. 
\label{stdfin}
\end{theorem}

{\em Note:} \hspace{.5mm}
Since the edges of
 our graph are geodesics in this case,
Morse theory assures us that we can assume that
there is an embedded disk $D$ in $B$ whose
interior is disjoint from
all of the edges of $\Gamma$ and whose boundary consists
of two arcs $\alpha$ and $\beta$, where $\alpha$ is
one of the edges of the graph, $\beta$ is strictly
contained in $\partial B$,
$\partial \alpha =
\partial \beta$ and the interior of 
$\beta$ is disjoint from all of the edges of $\Gamma$.

Let's examine how $D$ meets $\Sigma$.  

\begin{figure}[htbp] 
\centerline{\hbox{\psfig{figure=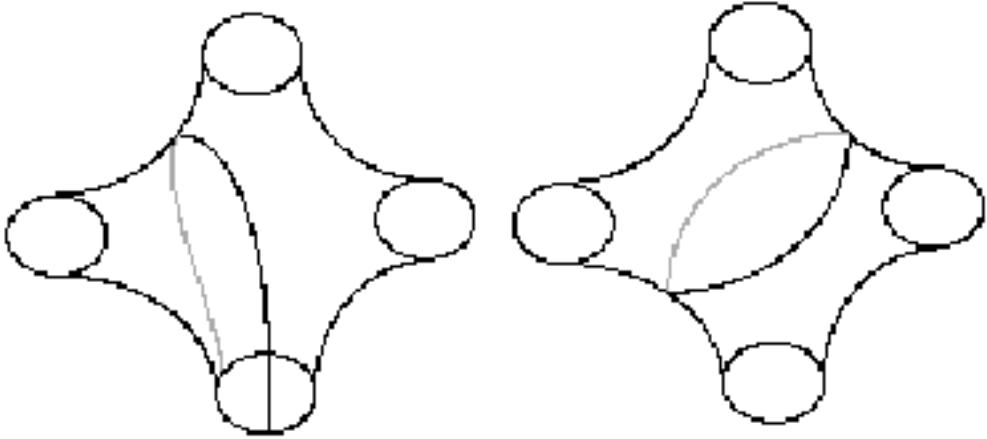}}}
\caption{An essential arc and curve
separate the remaining vertices into two sets.}\label{fig:separc}
\end{figure}

\begin{claim}
We may assume that
$D \cap \Sigma$ contains no simple closed curves.
\end{claim}
\begin{proof} Assume $D$ 
is chosen with a minimal number of intersections
with $\Sigma$.
Examine an innermost curve $\delta$ on $D$.  If $\delta$ is
not essential on $\Sigma$, then there is an
obvious isotopy through which this intersection
could have been eliminated, so we may assume it is essential on
$\Sigma$.   

Therefore the innermost loop gives us a compressing disk
for $\Sigma$.  Homology is enough to assure us that
the disk must be on the inside
of $\Sigma$ (the component of 
$B - \Sigma$ containing \ej).  Since $\delta$ is assumed to be
essential in $\Sigma$, it must separate the vertices
of $\Sigma$ into two non-empty sets.  
With no loss of generality, let $v_1$ be in one set and $v_2$
be in the other.
Since $\delta$ separates the vertices (as in
Figure~\ref{fig:separc}),  the disk it bounds does
too, and $e_{12}$ must intersect it.
This, however, is a contradiction since the interior of $D$
is disjoint from the edges of $\Gamma$.
\end{proof}

\vspace{.5in}

\begin{figure}[htbp] 
\centerline{\hbox{\psfig{figure=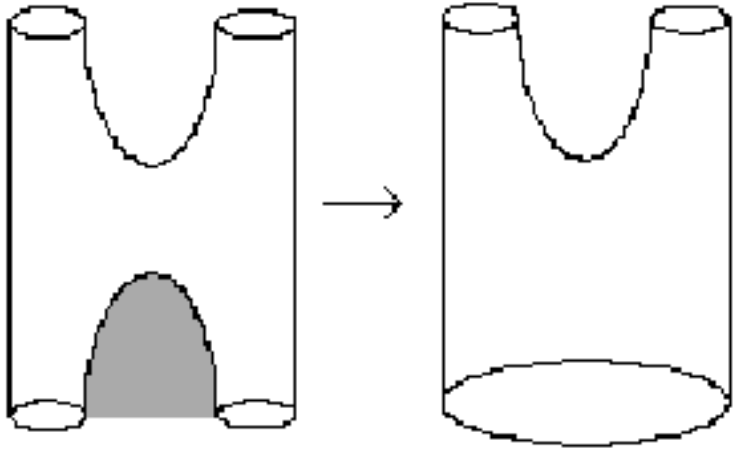}}}
\caption{A boundary-compressing disk on the outside
between two vertices reduces $n$ to $n-1$}
\label{fig:reduction}
\end{figure}

Now we examine an outermost arc $\gamma$ on $D$.
If $\gamma$ runs from one vertex of $\Sigma$
to a different one, then it is obvious that the
corresponding subdisk of $D$ is on the outside of $\Sigma$.
This gives us a compression disk that allows us to complete
the proof by induction. (See Figure~\ref{fig:reduction}).

We may therefore assume that $\gamma$ connects a vertex,
say $v_n$ to 
itself.  Because $D \cap \Sigma$ is assumed to be minimal
$\gamma$ must be essential on $\Sigma$.
If the disk $\gamma$ cuts off on $D$
 is on the inside of $\Sigma$,
then once again it must separate {$ \{ v_1, \dots v_{n-1} \}$  
into two sets and
the argument proceeds as in the simple closed curve case.
(See Figure~\ref{fig:separc}).

Our final case therefore is that although $\gamma$ connects
$v_n$ to itself, the disk is on the outside of $\Sigma$.
Compressing along this disk splits the n-punctured sphere
$\Sigma$ into two pieces $\Sigma_1$ and $\Sigma_2$.
$\Sigma_1$ is an r-punctured sphere and $\Sigma_2$
is a $n+1-r$-punctured sphere, where $ 2 \leq r \leq n-1$.
By induction we may therefore
assume that $\Sigma_1$ and $\Sigma_2$ are standard.
 $\Sigma_1$ and $\Sigma_2$ and their product structures
 are in different ``halves'' of $B$.
They are separated by the annulus formed by 
the boundary compressing $v_n$.
The inverse of the boundary compression is a tunnel
connecting the two boundary components of the annulus.
Up to isotopy there is a unique arc running across the annulus,
so there is a unique tunnel we can add to attain
$\Sigma$ from $\Sigma_1$ and $\Sigma_2$.
Since the product structures on $\Sigma_1$ and $\Sigma_2$
are to the outside of the annulus, and the tunnel can 
be added extremely close to the boundary, it is clear that
the product structure can be extended across the tunnel
to give a product
structure on $\Sigma$, proving that it is standard.

{\em Note:} \hspace{.5mm}we needed 
the full strength of the complete graph here.
Figure~\ref{fig:unstdsphere1} shows a counter-example for $n = 3$,
if one edge is missing from $\Gamma$.  
This counter-example may be generalized
for any $n$ to give a counterexample for the complete graph on 
$n$ vertices, minus one edge.

\begin{figure}[htbp] 
\centerline{\hbox{\psfig{figure=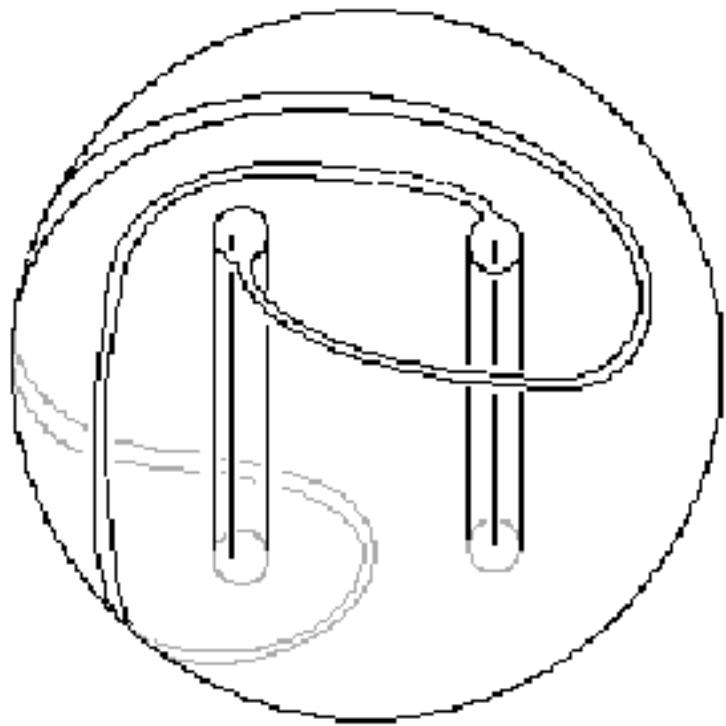}}}
\caption{A non-standard sphere for 
the complete graph on three vertices
minus one edge.}\label{fig:unstdsphere1}
\end{figure}

\section{The Infinite Case}

\begin{definition}
$f: F \rightarrow \Sigma \subset \R^3$ is 
a proper embedding of $F$, a sphere with a Cantor Set worth of 
punctures, in $\R^3$ if the pre-image of every compact set in $\R^3$ 
 is
a compact set in $F$ ($H^3$ could replace $\R^3$ throughout this section). 
\end{definition}

\begin{definition}
Let $\cup \it{l}_{\alpha}$ be a collection of geodesics in $\R^3$.
Let $\cup \it{l}_{\alpha}$ be contained in one
connected component of $\R^3 - \Sigma$.
Again let
the inside of $\Sigma$ be the component of $\R^3 - \Sigma$
containing $\Gamma$ and the outside be the remaining component. 
\end{definition}

\begin{definition} 
$\Sigma$ is said to be standard 
if there exists a product structure on the
outside of $\Sigma$ such that it is a product of the 
punctured sphere and a half open interval.

\end{definition}

We may now ask, how many lines contained inside
of $\Sigma$ it takes to ensure that the embedding of $\Sigma$ is standard.

\begin{theorem}  Given $\Sigma$, a proper embedding of $F$,
 a sphere with a Cantor Set worth of 
punctures in $\R^3$ and a set of geodesics which are dense in the Cantor Set,
we may conclude that the punctured sphere is standard.
\label{stdinf}
\end{theorem}

 In this context saying 
that the geodesics are dense in the Cantor Set
means that given any two punctures $p_1$
and $p_2$ of $F$ and any two neighborhoods
of those punctures $\mu_1$ and $\mu_2$ on $F$, there exists a geodesic 
that runs from the image of some puncture
in $\mu_1$ to the image of some puncture in $\mu_2$. 

Note that though there are an uncountable number of points in the 
cantor set, we are only requiring a countable set of geodesics.  
Even if we wanted to satisfy the property of being dense for every
point on the sphere, we would still only need a countable set of geodesics,
since geodesics connecting all the points on the sphere with rational
coordinates would suffice and the
set of possible pairings of a countable set
is itself a countable set.

Perhaps the easiest way to picture the scenario is to imagine
the universal cover of a genus-two handlebody in hyperbolic
three-space.  The punctured sphere would be the boundary of 
this cover and the lines would be geodesics connecting points
at infinity.  

It is worthy of note that in the finite case we
needed the complete graph, so we needed all possible geodesics,
but in the infinite case with an uncountable number of
punctures, we only need a countable number of edges.

\vspace{.25in}

We now begin the proof of the theorem.

\vspace{.25in}
\begin{proof} Choose a point on
$\Sigma$ to be the origin of $\R^3$.  

Let $S_n$ be the sphere of radius $n$ centered at the origin in
$\R^3$, and let $B_n$ be the ball that it bounds.
 We may alter $\Sigma$ slightly if necessary so
that we may assume that it intersects each $S_i$ transversally.

Fix $i$ and examine $B_i$.  The pieces of  $\Sigma$ in $B_i$
may be broken into two sets.  The first set consists of the
connected piece containing the 
origin $\Sigma_{i1}$, and the second set, all the other pieces, $\{\Sigma_{i2}, \dots
\Sigma_{in} \}$.  We shall isotope $\Sigma$
until $B_i$ contains one piece of the first type and none
of the second on the ``outside" of $\Sigma$ in $B_i$.  The latter
does not prevent us from claiming that the former is boundary parallel
in $B_i-\Sigma$, so we do not worry about them. See 
Figure~\ref{fig:cleanup}

\begin{figure}[htbp] 
\centerline{\hbox{\psfig{figure=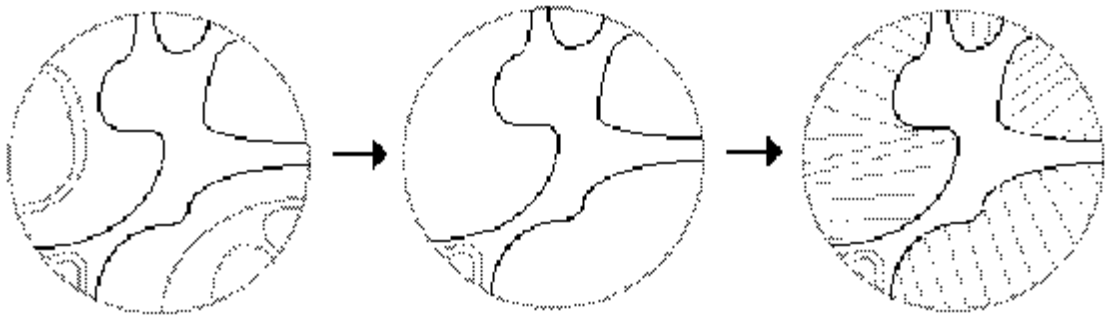,width=6in}}}
\caption{A small ball can be cleaned up leaving a product
structure on the outside of the punctured sphere.}\label{fig:cleanup}
\end{figure}

We examine the pre-images of the $\Sigma_{ij}$ in $F$.  
We now make an argument to show that we may assume that none of them
have a boundary curve which is trivial in the 
fundamental group of $F$.  

If there were a trivial boundary curve, we could choose an innermost
one (on $F$).  This would bound an embedded disk $D$ on $\Sigma$ that
meets $S_i$ in a simple closed curve.  The boundary curve splits
$S_i$ into two disks, each of which bounds a ball with $D$.
If $D$ is on the interior of $B_i$, then we choose the ball 
that does not contain the surface $\Sigma_{i1}$  ($\Sigma_{i1}$, is connected
and disjoint from $D$, so it can only be in one of the two balls).
If $D$ is on the exterior of $B_i$ then we choose the smaller
of the two
balls (it is contained in the other ball).  Either way we push $D$
across the ball and through $S_i$, eliminating at least its intersection
with $S_i$ and possibly more extraneous intersections that were contained
in the ball through which our isotopy was done.  

Note that $\Sigma_{i1}$ is unchanged away from its boundary and its 
boundary can only be changed by capping off trivial components.
We continue this process until there are no more trivial components
in the entire collection of $\Sigma_{ij}$.

At the risk of sloppy notation we shall continue 
throughout to
call the new surfaces $\Sigma_{ij}$ carefully noting 
at each step that we still have done only
a finite number of isotopies to a finite number of pieces.

We now notice that $\Sigma_{i1}$ fits all of the criterion of
the standard finite case.  Since the geodesics are dense within
the Cantor Set, at any finite stage there will be a complete graph in $B_i$
on the vertices that are given by the intersection of $\Sigma_{i1}$
and $S_i$.
Thus, by the previously proven finite case, 
$\Sigma_{i1}$ is standard in $B_{i}$.
We can use its inherited product structure to 
isotope $\Sigma$ to make sure that it
does not intersect $B_{i}$ on the outside of $\Sigma_{i1}$.

Now we repeat the process for some $k > i$.  We might worry that this 
process results in our pushing some piece of $\Sigma$ an infinite
number of times, but this is not the case, as every
point in $\Sigma$ is in some $\Sigma_{k1}$ for large enough $k$ and
our isotopies never affect points of $\Sigma_{k1}$ that are not
near the boundary of $B_k$.  Thus, for each point 
there is some $k$ such that the point is left alone for
good after $k$ steps.  Since each step involved only a finite
number of isotopies, each point is moved only a finite number of times.

The only thing left for us to check is that the product structure for
$\Sigma_{k1}$ can be chosen to correspond exactly with the product
structure we already chose for $\Sigma_{i1}$.  $B_i$
may be left fixed as we do our operations for $\Sigma_{k1}$,
so naturally $\Sigma_{i1}$ remains fixed, too.

Since $\Sigma_{i1}$ is boundary parallel in $B_i$, we may substitute
part of $S_i$ for it in $\Sigma_{k1}$ and the resulting surface
still contains a complete graph on one side and must be boundary
parallel.  If we concatenate its product lines in $B_k$ with the 
product lines of $\Sigma_{i1}$ in $B_i$ we see a product structure
that suits our desires as in Figure~\ref{fig:productextends}.

\begin{figure}[htbp] 
\centerline{\hbox{\psfig{figure=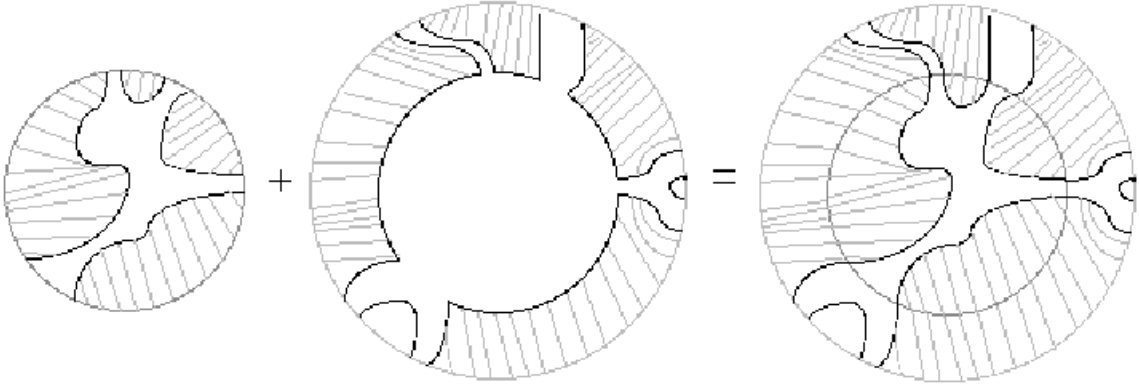,width=6in}}}
\caption{The product extends from one stage to the 
next}\label{fig:productextends}
\end{figure}
\end{proof}

\pagebreak
\section{References}
\hspace{.2in}

[B] J. Berge,
{\it The knots in $D^2 \times S^1$ which have nontrivial Dehn 
surgeries that yield $D^2 \times S^1$.} 
Topology Appl. 38 (1991), no. 1, 1--19.  

\vspace{.25in} 

[BZ] G. Burde and H. Zieschang,
{\it Knots}. Walter de Gruyter, Inc.  New York, NY, 1985.

\vspace{.25in} 

[C] J. H. Conway, {\it Enumeration of knots and links}. Computational
Problems in Abstract Algebra, Pergamon Press, 1970,  329--359.

\vspace{.25in}

[CF] J. Cannon and C. Feustel {\it Essential embeddings
of annuli and mobius bands in 3-manifolds}. Transactions
of the AMS 215, (1976), 219--239.

\vspace{.25in}

[GF] J. Goldman and L. Kauffman,
{\it Rational Tangles}. Advances in Applied Mathematics 18 (1997), 300--332.

\vspace{.25in}

[HR] J. Hempel and L. Roeling,
{\it Free factors of handlebody groups}.
preprint.

\vspace{.25in}

[M] K. Murasugi, {\it Knot Theory and Its Applications}, Birkhauser, Boston,
1996.
\vspace{.25in}

[ST] H. Schubert,
{\it Knoten mit zwei Brucken.} Math. Zeit.
1956, no. 65, 133--170

\vspace{.25in}

[ST] M. Scharlemann and A. Thompson,
{\it Detecting unknotted graphs in $3$-space.} J. Differential Geom. 
1991, no. 2, 539--560

\vspace{.25in}

[W] Jeff Weeks,
{\it SnapPea: a computer program for creating and studying 
hyperbolic 3-manifolds}, available by anonymous
        ftp from geom.umn.edu.

\end{document}